\documentclass[11pt]{article}
\usepackage{amsmath}
\usepackage{graphics, amssymb}
\newtheorem{Def}{Def.}[section]
\newtheorem{Thm}[Def]{Theorem}
\newtheorem{Prp}[Def]{Proposition}
\newtheorem{Lemma}[Def]{Lemma}

\newcommand{\Proof}{{\em{Proof. }}}
\newcommand{\QED}{\ \hfill $\FBox$ \\[1em]}
\newcommand{\spc}{\;\;\;\;\;\;\;\;\;\;}

\newcommand{\R}{\mathbb{R}}

\newcommand{\Tr}{\mbox{Tr\/}}

\newcommand{\beq}{\begin{equation}}
\newcommand{\eeq}{\end{equation}}

\newcommand{\FBox}{\rule{2mm}{2.25mm}}

\newcommand{\ov}{\overline}

\numberwithin{equation}{section}

\title{Curvature Estimates in Asymptotically Flat
Lorentzian Manifolds}
\author{Felix Finster, Margarita Kraus}
\date{June 2003}

\begin{document}
\maketitle

\begin{abstract}
We consider an asymptotically flat Lorentzian manifold of
dimension $(1,3)$. An inequality is derived which bounds the
Riemannian curvature tensor in terms of the ADM energy in the
general case with second fundamental form. The inequality
quantifies in which sense the Lorentzian manifold becomes flat in
the limit when the ADM energy tends to zero.
\end{abstract}

\section{Introduction}
In general relativity, space-time is modeled by a Lorentzian
manifold $(N,g)$ of signature $(- +++)$. Gravity is described
geometrically by Einstein's equations
\[ \overline{\mbox{Ric}} - \frac{1}{2}\: \overline{s}\: g \;=\;
-8 \pi\: T \;, \] where $\overline{R}$ is the curvature
corresponding to the Levi-Civita connection $\overline{\nabla}$ on
$N$, $\ov{\mbox{Ric}}$ is the Ricci curvature and $\ov{s}$ the
scalar curvature. Here $T$ is the energy-momentum tensor; it tells
about the distribution of matter in space-time and gives a local
concept of energy and momentum. The fact that the local energy
density should be positive is expressed by the {\em{dominant
energy condition}}, saying that for each $p \in N$ and each
timelike vector $u \in T_pN$, \begin{equation}
 T_{\alpha \beta}\: u^\beta
{\mbox{ is timelike}} \qquad {\mbox{and}} \qquad T(u,u) \leq 0\;.
\label{dec}\end{equation}
 We choose a space-like hypersurface $M
\subset N$ and let $(g,h)$ be the induced Riemannian metric and
the second fundamental form on $M$, respectively. In many physical
situations, matter is localized in a bounded region of space, and
the gravitational field falls off at large distance from the
sources. This leads to the definition of asymptotic flatness; for
simplicity we consider only one asymptotic end.
\begin{Def}
$M$ is {\bf{asymptotically flat}} if
there is a compact set $K \subset M$ and a diffeomorphism $\Phi$ which maps
$M \setminus K$ to the region $\R^3 \setminus B_r(0)$ outside a ball of radius $r$.
Under this diffeomorphism, the metric and second fundamental form should
be of the form
\begin{eqnarray*}
(\Phi_* g)_{ij} &=& \delta_{ij} + {\cal{O}}(r^{-1}) \;,\quad
\partial_k (\Phi_* g)_{ij} \;=\; {\cal{O}}(r^{-2}) \;,\quad
\partial_{k}\partial_l (\Phi_* g)_{ij} \;=\; {\cal{O}}(r^{-3}) \\
(\Phi_* h)_{ij} &=& {\cal{O}}(r^{-2}) \;,\qquad\;\;\;\;\,
\partial_k (\Phi_* h)_{ij} \;=\; {\cal{O}}(r^{-3}) \;.
\end{eqnarray*}
\end{Def}
In asymptotically flat manifolds, one can introduce the ADM energy
and momentum, which have the interpretation as the total energy and momentum
of space-time.
\begin{Def} The {\bf{ADM energy and momentum}} $(E,P)$ are defined by
\begin{eqnarray}
E &=& \frac{1}{16 \pi} \lim_{R \to \infty} \sum_{i,j=1}^3
\int_{S_R} (\partial_j (\Phi_* g)_{ij} - \partial_i (\Phi_*
g)_{jj}) \:d\Omega^i \label{12} \\ P_k &=& \frac{1}{8 \pi} \lim_{R
\to \infty} \sum_{i=1}^3\int_{S_R} ((\Phi_* h)_{ki} - \sum_{j=1}^3
\delta_{ki} \:(\Phi_* h)_{jj}) \:d\Omega^i \;\;\; , \label{13}
\end{eqnarray}
where $d\Omega^i = \nu^i \:du$, $du$ is the area form, and $\nu$ is
the normal vector to $S_R \subset \R^3$.
\end{Def}
This definition is indeed independent of the the choice of
$\Phi$~\cite{ADM}.

It is a major problem of mathematical relativity to understand the
relation between $(E, P)$ and the geometry of space-time. A
particular aspect of this problem is the question whether and in
which sense $E$ and $P$ control the Riemannian curvature tensor.
In~\cite{BF} this question was addressed in the time-symmetric
case (i.e.\ when $h \equiv 0$). $L^2$-estimates for the Riemannian
curvature tensor where derived on $M \setminus D$, where $D$ is an
``exceptional set'' of small volume. In~\cite{FK} these estimates
were generalized to higher dimensions. In the present paper we
treat the physically relevant case with second fundamental form.
This is our main result:
\begin{Thm}
We choose $L \geq 3$ such that
\[ (L^\alpha -1)^2 \;\geq\; C\: \frac{4\pi E+\|h\|_2}{k^2\: (k+24\: \|h\|_3)^2} \: \||h|^2 + |\nabla h|\|_3 \]
where
\[ \alpha \;=\; \left( 1 + 24\: \frac{\|h\|_3}{k} \right)^{-1} . \;\]
Then there is a set $U$ with measure bounded by
\[ \mu(U) \;\leq\; c_1\: \frac{L^6}{k^2} \:(4\pi E + \|h\|^2_2) \]
such that on $M \setminus U$ the following inequality holds,
\begin{eqnarray*}
\lefteqn{ \int_{M \setminus U} \eta\: |\overline{R}_M|^2\:d\mu \;\leq\;
c_2\: \sup_{M}\left( |\Delta\eta|+|\nabla\eta||h|+\eta\: (|R|+|h|^2+|\ov{\nabla}h| \right) E } \\
&&+c_3\:L \:\sup_{M} \left( \eta\:(|\ov{\nabla}\ov{R}_M|+|h||\ov{R}_M|)\right) \sqrt{E} \\
&&+c_4\:\frac{\sqrt{L+1}}{k}
\:(\sup_{M}\eta)\; \sqrt{\||h|^2+|\nabla h|\|_{6/5}} \; \left\||\ov{\nabla}\ov{R}_M|+|h||\ov{R}_M \right\|_{5/12}
\;\sqrt{E}\;.
\end{eqnarray*}
Here $c_1,\ldots,c_4$ are numerical constants (independent of $L$ and the geometry),
$\eta \in C^2(M)$ is a positive test function, $\overline{R}_M$ is the
Riemannian curvature tensor of $N$ restricted to $M$, and $k$ is the
isoperimetric constant $k=\inf A/V^{\frac{2}{3}}$.
\end{Thm}

For the proof we use Witten's solutions of the hypersurface Dirac
equation~\cite{W, PT} and consider second derivatives of the
spinors. In order to control the Weyl tensor, we work similar as
in~\cite{FK} with the spinor operator $\Pi$, which is built up of
a whole family of solutions of the hypersurface Dirac equation.
The presence of the second fundamental form leads to the
difficulty that the function $|\Psi|^2$ is no longer subharmonic,
making it impossible to estimate the norm of the spinor with the
maximum principle. In order get around this difficulty, we first
construct a barrier function $F$, which is a solution of a
suitable Poisson equation. We then derive Sobolev estimates for
$F$, and these finally give us control of
$\||\Psi|^2-1\|_{L^6(M)}$.

\section{Basic Facts about Spinors and the Hypersurface Dirac Operator}
In this section we recall some basic facts about spinors and the
Dirac operator on Lorentzian spin manifolds; for details the
reader is referred to \cite{Baum}, \cite{Lami}.

Let $(N,g)$ be a Lorentzian spin-manifold with spin structure
$Q_N\to N$. Let\linebreak $\kappa
:\mbox{Spin}(1,n-1)\to\Delta_{1,n-1}$ denote the spinor
representation and
\begin{equation*}\Sigma_N \;= \;Q_N\times_\kappa\Delta_{1,n-1}
\end{equation*}
the associated spinor bundle.
 We denote the Clifford
multiplication of a tangent vector $X$ with a spinor $\psi$ by
$\mu(X,\psi)=:X\cdot\psi$. On $\Delta_{1,n-1}$ there exists an
indefinite scalar product $\langle\cdot,\cdot\rangle$ of signature
$(2,2)$, which is invariant under $\mbox{Spin}^+(1,n-1)$ and is
unique up to a constant. This inner product induces on $\Sigma_N$
an indefinite scalar product, which we again denote by
$\langle\cdot,\cdot\rangle$. For a timelike vector field $\nu$,
the inner product
\begin{equation}
(\phi,\psi) := \langle\phi,\nu\cdot \psi\rangle \label{sp}
\end{equation}
is positive.

The scalar products $(\cdot,\cdot)$ and
$\langle\cdot,\cdot\rangle$ also define scalar products on the
fibres of the bundles of $k$-forms $\Lambda^k T^\ast N\otimes
\Sigma_N$ and the bundle of $k$-linear mappings $(\otimes^kT^\ast
N\otimes\Sigma_N)$ by
\begin{equation*}
(\eta,\xi)_p \; :=\;
\sum_{i_1,\dots,i_k}(\eta(e_{i_1},\dots,e_{i_k}),
\xi(e_{i_1},\dots,e_{i_k}))_p
\end{equation*}
and analogously
\begin{equation*}
\langle \eta,\xi\rangle_p \;:=
\;\sum_{i_1,\dots,i_k}\langle\eta(e_{i_1},\dots,e_{i_k}),\xi(e_{i_1},\dots,e_{i_k})\rangle\;,
\end{equation*}
where $e_1,\dots,e_n$ is an orthonormal frame. The Levi-Civita
connection on $N$ induces a covariant derivative $\ov{\nabla}$ on
$\Gamma\Sigma_N$.

This covariant derivative is isometric with respect to
$\langle\dots\rangle$, i.e.
\begin{equation*}
X\langle\varphi,\psi\rangle \; =\;
\langle\ov{\nabla}_X\varphi,\psi\rangle+\langle\varphi,\ov{\nabla}_X\psi\rangle
\end{equation*}
for all sections $\varphi,\psi$ in $\Sigma_N$. Its curvature
tensor $\ov{R}^\Sigma\in\Omega^2(N,\Sigma_N)$ is defined by
\begin{equation*}
\ov{R}^\Sigma (X,Y)\psi \; =\;
(\ov{\nabla}^2\psi)(X,Y)-(\ov{\nabla}^2\psi)(Y,X),
\end{equation*}
where
$(\ov{\nabla}^2\psi)(X,Y)=\ov{\nabla}_X\ov{\nabla}_Y-\ov{\nabla}_{\ov{\nabla}_XY}$.
It is related to the curvature tensor $\ov{R}$ of the Lorentzian
manifold $(N,g)$ by the formula
\begin{equation}
\ov{R}^\Sigma \psi\; =
\;\frac{1}{4}\sum^n_{\alpha,\beta=1}\langle\ov{R}e_\alpha,e_\beta\rangle
e_\alpha\cdot e_\beta\cdot\psi .\label{curv}
\end{equation}
The Dirac operator on the Lorentzian manifold $N$ is defined by
the composition of the covariant derivative $\ov{\nabla}$ with the
Clifford multiplication $\mu$,
\begin{equation*}
\ov{D}:\Gamma\Sigma_N\stackrel{\ov{\nabla}}{\longrightarrow}\Gamma(T^\ast
N\otimes\Sigma_N)\stackrel{\mu}{\longrightarrow}\Gamma\Sigma_N,
\end{equation*}
where the cotangent bundle $T^\ast N$ has been identified with the
tangent bundle $TN$ via the metric. In a local orthonormal frame
$(e_1,\dots,e_n)$, the Dirac operator is given by
\begin{equation*}
\ov{D} \; =
\;\sum^n_{\alpha=1}e_\alpha\cdot\ov{\nabla}_{e_\alpha}.
\end{equation*}
We point out that, in contrast to the Riemannian case, the Dirac
operator on a Lorentzian  manifold is not elliptic.

In what follows, we restrict attention to the physically relevant
case of a $4$-dimensional Lorentzian manifold with a given
$3$-dimensional asymptotically flat space like hypersurface
$M\subset N$. We choose a normal unit vector field $\nu$ on $M$
and consider the corresponding positive definite scalar pro\-duct
(\ref{sp}). We set $\vert\psi\vert = (\psi,\psi)^{\frac 12}.$

The covariant derivative $\ov{\nabla}$ is not compatible with
$(\cdot,\cdot)$, but
\begin{equation*}
X(\varphi,\psi)\; =
\;(\ov{\nabla}_X\varphi,\psi)+(\varphi,\ov{\nabla}_X\psi)+(\varphi,\nu\cdot\ov{\nabla}_X\nu\cdot\psi)\label{lr}
\end{equation*}
holds for spinor fields $\varphi,\psi\in\Gamma(\Sigma_N|M)$. Using
the definition of the second fundamental form
\begin{equation*}
h_{ij} \; = \;-\langle e_i,\ov{\nabla}_{e_j}\nu\rangle
\end{equation*}
for an orthonormal frame $(e_1,e_2,e_3)$ on $M$, this formula can
be written as
\begin{equation}\label{der}
e_i(\varphi,\psi) \;=
\;(\ov{\nabla}_{e_i}\varphi,\psi)+(\varphi,\ov{\nabla}_{e_i}\psi)-h_{ij}(\varphi,\nu\cdot
e_j\cdot\psi).
\end{equation}
This leads us to define the adjoint of $\ov{\nabla}$ by
$\ov{\nabla}_X^\ast\psi=-\ov{\nabla}_X\psi-\nu\cdot\ov{\nabla}_X\nu\cdot\psi$
or, in an orthonormal frame,
\begin{equation*}
\ov{\nabla}_{e_i}^\ast\psi\; =
\;-\ov{\nabla}_{e_i}\psi+h_{ij}\nu\cdot e_j\cdot\psi.
\end{equation*}
On a spacelike hypersurface, there exists an intrinsic Riemannian
Dirac operator, but we shall not consider it here. Instead, we
will only be concerned with the so called {\bf hypersurface Dirac
operator} $\overline{D}_M$,
\begin{equation*}
\ov{D}_M\;:= \;\ov{D}\vert_M:\Gamma
\Sigma_{\ov{M}}|M\to\Gamma\Sigma_{\ov{M}}|M.
\end{equation*}
It is the restriction of the Dirac operator of the Lorentzian
manifold $N$ to $M$; more precisely,
\begin{equation*}
\Gamma(\Sigma_{\ov{M}}\vert
M)\stackrel{\ov{\nabla}}{\longrightarrow} \Gamma(T^\ast
M\otimes\Sigma_{\ov{M}}|M)\longrightarrow
\Gamma(\Sigma_{\ov{M}}|M),
\end{equation*}
where $\ov{\nabla}$ denotes the covariant derivative in direction
$M$. According to \cite{PT}, the square of the hypersurface Dirac
operator satisfies the Weitzenb{\"o}ck formula
\begin{equation}
\ov{D}^2_M\; =\;
\ov{\Delta}^s+\Re. \label{wb}
\end{equation}
Here $\ov{\Delta}^s$ is the Laplacian $\ov{\Delta}^s
\psi=\ov{\nabla}^\ast\ov{\nabla}\psi=
\mbox{tr}(-\ov{\nabla}^2-\nu\cdot\ov{\nabla}\nu\cdot\ov{\nabla}\psi)$
or, in an orthonormal frame,
\begin{equation*}
\ov{\Delta}^s\psi \; =
\;-\sum_{i,j}(\ov{\nabla}_{e_j}\ov{\nabla}_{e_j}\psi-\ov{\nabla}_{\nabla_{e_j}e_j}\psi
-  h_{ij}\nu\cdot e_i\cdot\ov{\nabla}_{e_j}\psi),
\end{equation*}
and $\Re$ is the curvature expression $\Re
=\frac{1}{4}(\ov{s}+2\,\ov{\mbox{Ric}}(\nu,\nu)+2\sum^3_{i=1}\ov{\mbox{Ric}}(\nu,e_i)(\nu\cdot
e_i)$. The dominant energy condition~(\ref{dec}) yields that $\Re
\geq 0$.

In the coordinates induced by the diffeomorphism $\phi$ of
Definition 1.1, we choose a constant spinor $\psi_0$ of norm one
in the asymptotic end and consider the boundary value problem
\begin{equation}\label{br}
\ov{D}_M\psi \; =\; 0\,,\qquad \lim_{|x|\to\infty}\psi(x)\; =
\;\psi_0 {\mbox{ with }} |\psi_0|=1\:.
\end{equation}
The existence and uniqueness of a solution of (\ref{br}) is proven
in \cite{PT}. The solution decays at infinity as
\begin{equation*}
\psi\;=\;\psi_0+\mathcal{O}(r^{-1}),\enskip
\partial_j\psi\;=\;\mathcal{O}(r^{-2}),\enskip
\partial_{kl}\psi\;=\;\mathcal{O}(r^{-3}).
\end{equation*}
Using the Weitzenb\"{o}ck formula~(\ref{wb}),  it is shown in
\cite{PT} that for a solution of (\ref{br}),
\begin{equation} \label{abs}
\|\nabla\psi\|^2_{L^2(M)}\;=\;4\pi \left(E\:|\psi_0|^2+(\psi_0,\:P\cdot\psi_0) \right) -(\psi,\:\Re\psi)
\;\leq\; 4\pi \left( E+(\psi_0,\:P\cdot\psi_0) \right)\;,
\end{equation}
where $P =P_k\cdot e_k$ is the momentum as defined by (1.2). If we
choose $\psi_0$ such that $(\psi_0|P\cdot\psi_0)=-|P|$, we obtain
the positive mass theorem~\cite{W, PT}
\begin{equation} \label{pmt}
0 \;\leq\; 4\pi\:(E-|P|)\;.
\end{equation}
For general $\psi_0$, (\ref{abs}) and~(\ref{pmt}) give rise to
an $L^2$-bound of $\nabla \psi$,
\begin{equation} \label{mf}
\|\nabla\psi\|^2_{L^2(M)} \;\leq\; 4 \pi\: (E+|P|) \;\leq\; 8 \pi E \;.
\end{equation}

\section{A-priori Estimates for Harmonic Spinors}

In what follows, we let $\psi\in\Gamma \Sigma_N|M$ be a solution
of the boundary value problem (\ref{br}). We refer to $\psi$ as a
{\em{harmonic spinor}}. We begin by deriving an upper bound for
the measure of the set where a harmonic spinor is large. For any
$L \geq 1$, we introduce the set $\Omega_L = \Omega_L (\psi)$ by
\begin{equation}
 \Omega_L (\psi) \;=\; \{x\in
M \,:\, |\psi(x)|\geq L\} \;.\label{dol}
\end{equation}

\begin{Lemma}\label{volume}
For any harmonic spinor $\psi$ and every $L \geq 1$, the volume of
$\Omega_L$ is bounded by
\[ \mu(\Omega_L)^{\frac{1}{3}} \;\leq\; \frac{192}{(L^\alpha-1)^2}\;
\frac{4 \pi E + \|h\|_2^2}{k^2}\;, \] where the exponent $\alpha$
is
\begin{equation} \label{alpha}
\alpha \;=\; \left( 1 + 24\: \frac{\|h\|_3}{k} \right)^{-1} .
\end{equation}
\end{Lemma}
The proof uses the the following Sobolev inequality, which is
derived in~\cite{FK}.
\begin{Lemma} \label{lemma31}
Let $M$ be an asymptotically flat manifold of dimension $n \geq 3$.
Then every non-negative function $g\in C^\infty(M) \cap H^{1,2}(M)$ with
$\underset{|x|\to\infty}{\lim} g(x)=0$ satisfies the inequality
\[ \|g\|_q \;\leq\; \frac{q}{k}\: \|\nabla g\|_2 \spc {\mbox{with}} \spc
q = \frac{2n}{n-2}  \] and $k$ the isoperimetric constant.
\end{Lemma}
{\em{Proof of Lemma~\ref{volume}.}} Applying the Schwarz
inequality in~(\ref{der}), we obtain for every $\alpha \in \R$,
\[ |\nabla|\psi|^{\alpha}| \;\leq\;
\alpha|\psi|^{\alpha-2}(|\nabla\psi|\:|\psi|+|h|\:|\psi|^2) \;. \]
We take the square and use the inequality $(x+y)^2 \leq 2(x^2+y^2)$,
\[ |\nabla|\psi|^{\alpha}|^2 \;\leq\;
2\alpha^2 \left(|\nabla\psi|^2\: |\psi|^{2\alpha-2} \:+\:|h|^2\:
|\psi|^{2\alpha} \right) \;. \] Choosing $x\in\Omega_L$ and
$\alpha \in (0, 1]$, the factor $|\psi(x)|^{2\alpha-2}<1$, and
thus at $x$,
\[ |\nabla|\psi|^{\alpha}|^2 \;\leq\; 2 \alpha^2
\left(|\nabla\psi|^2 + |h|^2\: |\psi|^{2\alpha} \right) \;. \]
We integrate over $\Omega_L$ and apply Lemma~\ref{lemma31} as well as~(\ref{mf}),
\[ \||\psi|^{\alpha}-1\|_{L^6(\Omega_L)}^2 \;\leq\; \frac{72\:\alpha^2}{k^2}\:
\left(8 \pi E + \||h|^2\:|\psi|^{2\alpha}\|_{L^1(\Omega_L)}  \right) \; . \]

The last inequality has the disadvantage that the spinor also appears on the
right. Therefore, we apply the inequality
$|\psi|^{2\alpha}\leq 2(|\psi|^{\alpha}-1)^2+2$ and H{{\"o}}lder to obtain
\[ \||\psi|^{\alpha}-1\|_{L^6(\Omega_L)}^2 \;\leq\; \frac{72\:\alpha^2}{k^2}\:
\left(8 \pi E + 2\:\|h\|^2_3 \;\|(|\psi|^{\alpha}-1)\|^2_{L^6(\Omega^1)}
+ 2\:\|h\|^2_2 \right) \; . \]
Now we can combine the terms involving the spinors,
\begin{equation} \label{11}
\left[1 - \frac{12^2\:\alpha^2}{k^2}\:\|h\|^2_3 \right] \||\psi|^{\alpha}-1\|_{L^6(\Omega_L)}^2
\;\leq\; \frac{12^2\:\alpha^2}{k^2}\:
\left(4 \pi E + \|h\|^2_2 \right) \; .
\end{equation}

We choose $\alpha$ according to~(\ref{alpha}).
Then the second term in the square brackets in~(\ref{11}) is bounded by
\[ \frac{12^2\:\alpha^2}{k^2}\:\|h\|^2_3 \;\leq\; \frac{1}{4} \]
and thus
\[ \||\psi|^{\alpha}-1\|_{L^6(\Omega_L)}^2 \;\leq\; 192\: \frac{\alpha^2}{k^2}\: \left(4 \pi E + \|h\|^2_2 \right) \;. \]
We finally apply the estimate
\[ \mu(\Omega_L)^{\frac{1}{3}} \;\leq\; \frac{1}{(L^\alpha-1)^2}\: \||\psi|^\alpha - 1\|_{L^6(\Omega_L)}^2\;. \]
\QED

In the time-symmetric case, Lemma~\ref{volume} reduces to the
inequality
\begin{equation} \label{1es}
\mu(\Omega_L)^{\frac{1}{3}} \;\leq\; \frac{192}{(L-1)^2}\; \frac{4 \pi E}{k^2}\;,
\end{equation}
showing that for large $L$, $\mu(\Omega_L)$ decays at least $\sim
L^{-6}$. On the other hand, it was shown in the time-symmetric
case~\cite{BF} that the function $|\psi|^2$ is subharmonic, and
thus the maximum principle gave the bound
\begin{equation}
|\psi|^2 \;\leq\; 1\;. \label{1et}
\end{equation}
This shows that if $h \equiv 0$, $\mu(\Omega_L)$ is indeed zero
for all $L>1$. We conclude that the estimate~(\ref{1es}) is
certainly not optimal if $h \equiv 0$. We shall now improve
Lemma~\ref{alpha} such that in the time-symmetric case we
recover~(\ref{1et}). We let $(e_1,e_2,e_3)$ be an orthonormal
frame in a neighborhood of $x$ with $(\ov{\nabla}_ie_j)(x)=0$.
Then the Laplacian of $|\psi|^2$ at $x$ is computed as follows,
\begin{eqnarray*}
 \Delta|\psi|^2 &=& \sum\limits^3_{j=1}
\partial_j((\ov{\nabla}_j\psi,\psi)+(\psi,\ov{\nabla}_j\psi)+(\psi,\nu\cdot\ov{\nabla}_j\nu\cdot\psi))
 \\
&=&2\:|\ov{\nabla}\psi|^2+2\:\mbox{Re}(\ov{\nabla}_j\ov{\nabla}_j\psi,\psi)+2\:(\ov{\nabla}_j\psi,\nu\cdot\ov{\nabla}_j\nu\cdot\psi)
\\ && +(\psi,\nu\cdot\ov{\nabla}_j\nu\cdot\ov{\nabla}_j\psi)
+(\psi,\ov{\nabla}_j(\nu\cdot\ov{\nabla}_j\nu\cdot\psi))-|\ov{\nabla}_j\nu|^2|\psi|^2
\\
&=&2|\ov{\nabla}\psi|^2-2\mbox{Re}(\ov{\nabla}^\ast\ov{\nabla}\psi,\psi)
+2\mbox{Re}(\ov{\nabla}_j\psi,\nu\cdot\ov{\nabla}_j\nu\cdot\psi)
+(\psi,\nu\cdot\ov{\nabla}^2_{j,j}\nu\cdot\psi))
\end{eqnarray*}
Using the Weitzenb{{\"o}}ck formula, we obtain for a harmonic
spinor the inequality $$|\Delta|\psi|^2| \;\geq\;
2\:\mbox{Re}(\Re\psi,\psi)+2|\ov{\nabla}\psi|^2-
2|\ov{\nabla}\psi||\psi||\ov{\nabla}\nu|-|\psi|^2\cdot\left(\sum\limits^3_{j=1}|
\ov{\nabla}^2_{j,j}\nu|\right) , $$ where we set
\[ |\ov{\nabla}\nu|^2 \;=\; \sum\limits_{i,j}h^2_{ij} \spc {\mbox{and}} \spc
|\ov{\nabla}^2_{j,j}\nu| \;=\;
3\sqrt{\sum\limits_k(\partial^jh_{jk})^2} \;. \] Using the short
notation
\[ |h| \;:=\; |\ov{\nabla}\nu| \spc {\mbox{and}} \spc
|\nabla h| \;:=\; \sum_{j=1}^3 |\ov{\nabla}^2_{j,j}\nu| \;,
\] we can write the last inequality in the compact form
\begin{eqnarray*}
\Delta|\psi|^2 &\geq&
2\mbox{Re}(\Re\psi,\psi)+2|\ov{\nabla}\psi|^2-2|\ov{\nabla}\psi||\psi||h|
-|\psi|^2|\nabla h|\\ &\geq&
2\mbox{Re}(\Re\psi,\psi)-\left(\frac{1}{2}|h|^2+|\nabla h|\right)
|\psi|^2\;.
\end{eqnarray*}
In the special case $h\equiv0$, we recover that $|\psi|^2$ is
subharmonic, and the maximum principle gives~(\ref{1et}). Our
method for treating the general case is to construct a barrier
function $F$ by solving the Poisson equation and to estimate $F$
using Sobolev techniques and the volume bound of
Lemma~\ref{volume}.
\begin{Prp} Suppose that $L>1$ is chosen so large that
\begin{equation} \label{Lbed}
(L^\alpha -1)^2 \;\geq\; C\: \frac{4\pi E+\|h\|_2}{k^2\: (k+24\: \|h\|_3)^2} \: \||h|^2 + |\nabla h|\|_3
\end{equation}
with $\alpha$ as in Lemma~\ref{volume} and $C=6\cdot 48^2$ a
numerical constant. Then the harmonic spinor $\psi$  is bounded on
$\Omega_L$ by
\[  \||\psi|^2-1\|_{L^6(\Omega_L)} \;\leq\; \frac{72}{k^2}\: (L+1)\:\||h|^2 + |\nabla h|\|_{L^{6/5}} \;. \]
\end{Prp}
{\Proof} We set $\rho=-(|h|^2+|\nabla h|)$ and let $g$ be the
solution of the Poisson equation $\Delta g=\rho|\psi|^2$ with
boundary conditions $\underset{x\to\infty}{\lim} g(x)=0$ (For the
existence of this solution see \cite[Theorem1.7]{Ba}). Then
$\Delta(|\psi|^2-g)\geq 0$, and the maximum principle yields that
\begin{equation} \label{1a}
|\psi|^2 \;\leq\; 1+g\,.
\end{equation}
The Sobolev inequality of Lemma~\ref{lemma31}, Gauss' theorem, and the
H{{\"o}}lder inequality give
\[ \|g\|^2_6 \;\leq\; \frac{36}{k^2}\|\nabla g\|^2_2
\;=\; \frac{36}{k^2}\int\limits_M |\rho||\psi|^2 \:g\:dM
\;\leq\; \frac{36}{k^2}\: \|\rho|\psi|^2\|_{\frac{6}{5}} \: \|g\|_6 \]
and thus
\begin{equation} \label{1b}
\|g\|_6 \;\leq\; \frac{36}{k^2} \:\|\rho|\psi|^2\|_{\frac{6}{5}}
\;.
\end{equation}
Combining~(\ref{1a}) and~(\ref{1b}), we obtain for any $L \geq 1$,
\begin{eqnarray*}
\||\psi|^2-1\|_{L^6(\Omega_L)} &\leq& \|g\|_{L^6(\Omega_L)}
\;\leq\; \frac{36}{k^2}\: \|\rho|\psi|^2\|_{\frac{6}{5}} \\ &\leq&
\frac{36}{k^2} \left(L\:\|\rho\|_{L^{6/5}(M\setminus\Omega_L)}
+\|\rho|\psi|^2\|_{L^{6/5}(\Omega_L)}\right)\\ &\leq&
\frac{36}{k^2} \left(L\:\|\rho\|_{L^{6/5}(M\setminus\Omega_L)}+
\|\rho\:(|\psi|^2-1)\|_{L^{6/5}(\Omega_L)} +
\|\rho\|_{L^{6/5}(\Omega_L)}\right)\\ &\leq& \frac{36}{k^2}
\left((L+1)\:\|\rho\|_{\frac{6}{5}}+
\|\rho\|_{L^{3/2}(\Omega_L)}\:
\||\psi|^2-1\|_{L^6(\Omega_L)}\right).
\end{eqnarray*}
We collect all the terms which involve $\||\psi|^2-1\|_{L^6(\Omega_c)}$,
\[ \left(1 - \frac{36}{k^2}\:\|\rho\|_{L^{3/2}(\Omega_L)} \right) \||\psi|^2-1\|_{L^6(\Omega_L)}
\;\leq\; \frac{36}{k^2} (L+1)\:\|\rho\|_{\frac{6}{5}} \;. \] This
inequality gives a bound for $\||\psi|^2-1\|_{L^6(\Omega_L)}$ only
if the prefactor is bounded away from zero. Thus we want to
arrange that
\begin{equation} \label{notw}
\frac{36}{k^2}\:\|\rho\|_{L^{3/2}(\Omega_L)} \;\leq\; \frac{1}{2}\;.
\end{equation}
The H{\"o}lder inequality gives
\[ \|\rho\|_{L^{3/2}(\Omega_L)} \;\leq\; \|\rho\|_3\: \mu(\Omega_L)^{\frac{1}{3}} \;. \]
Substituting in the volume bound of Lemma~\ref{volume}, one sees
that~(\ref{Lbed}) indeed guarantees that~(\ref{notw}) holds.
\QED

\section{Estimates of the Spinor Operator}

We choose an orthonormal basis of constant spinors
$(\psi^i_0)_{i=1,\dots,4}$, $(\psi^i_0,\psi^j_0)\equiv\delta_{ij}$
at the asymptotic end and denote the corresponding solutions of
the boundary problems~(\ref{br}) by $(\psi^i)_{i=1,\dots,4}$.

For every $x \in M$ we introduce the {\bf spinor operator} $\Pi_x$
by $$\Pi_x:\Sigma_{N,x} \to \Sigma_{N,x},
\psi\mapsto\sum\limits^4_{i=1} (\psi^i(x),\psi)\psi^i(x).$$ At
infinity, $(\psi^i)$ goes over to an orthonormal basis, and thus
$$\underset{|x|\to\infty}{\lim} \Pi_x=\mbox{id}$$ The next
elementary lemma bounds the spinor operator in terms of the
$|\psi^i|$.

\begin{Lemma} The $\sup$-norm of $\Pi_x$ is bounded by
\[ \frac{1}{4}\: \sum_{j=1}^4 |\psi^i_x|^2 \;\leq\; |\Pi_x| \;\leq\; \sum_{j=1}^4 |\psi^i_x|^2 \;. \]
\end{Lemma}
\begin{Proof} Since $\Pi_x$ is positive,
\[ |\Pi_x| \;\geq\; \frac{1}{4}\: {\mbox{tr}} \:\Pi_x \;=\; \frac{1}{4}\: \sum_{j=1}^4 |\psi^i_x|^2 \;. \]
This is the lower bound.

In order to derive the upper bound, we define the matrix $A$ by
$A=(a_{ij})_{\substack{i=1,\dots,4\\j=1,\dots,4}}$ with
$$a_{ij}=(\psi^i_x,\psi^j_x)\,.$$ By definition, $A$ is Hermitian
and all eigenvalues of $A$ are real and nonnegative. Let
$v:=(v_1,\dots,v_4)^T\in\mathbb{C}^4$, $|v|^2=1$. Then
$\psi:=\Sigma v_i\psi^i$ is a solution of the boundary problem
(\ref{br}) with $\psi_0=\Sigma v_i\psi^i_0$. Then $$
\sum\limits^4_{i,j=1}\ov{v}_iv_ja_{ij} \;=\; |\psi_x|^2 \;\leq\;
\sum_{j=1}^4 |\psi^i_x|^2 =: \lambda \;. $$ Therefore the
eigenvalues of $A$ must be smaller or equal  to $\lambda$.

Now let $\phi$ be an arbitrary spinor at $x \in M$. We let $\psi:=\Sigma v_i\psi^i_x$ be the
orthonormal projection of $\phi$ onto the span of
$(\psi^1_x,\dots\psi^4_x)$ and set $v^T=(v_1,\dots,v_4)$. Then
$$\begin{array}{lcl}
|\Pi_x\phi|^2&=&\sum\limits^4_{i,j=1}(\phi,\psi^i_x)
(\psi^i_x,\psi^j_x)(\psi^j_x,\phi)=(\ov{Av})^TA(Av) \\ &\leq &
\lambda^2\phantom{a}^T\ov{v}Av \;=\; \lambda^2\:|\psi_x|^2 \;\leq\; \lambda^2\:|\phi|^2
\end{array}
$$ and thus $|\Pi_x\phi|\leq \lambda |\phi|$.\QED
\end{Proof}

Next we derive an estimate for the Hilbert-Schmidt Norm
$\|\cdot\|$ of the operator $\|1-\Pi_x\|$.
\begin{Lemma} \label{lemma42}
For every $L \geq 3$ and $\varepsilon \in (0,1)$ there is a subset
$U \subset M$ with
\begin{equation*}
\mu (U)^{\frac{1}{3}} \leq \frac{48}{k^2} (4\pi E +
\|h\|^2_2)\;\frac{L^2(4+L^2)^2}{\varepsilon^2}
\end{equation*}
such that for all $x \in M \setminus U,$
$$\|1-\Pi_x\|<\varepsilon.$$
\end{Lemma}
\begin{Proof} We set $p(x) =\|1-\Pi_x\|^2$. Then the same calculation as in
  \cite[Lemma 4.2]{FK}, shows that
$$p(x)=4-2\sum\limits^4_{i=1}(\psi^i_x,\psi^i_x)+\sum\limits^4_{i,j=1}|(\psi^i_x,\psi^j_x)|^2.$$
Differentiation gives $$\begin{array}{lcl} \nabla p&=&
-4\sum\limits^4_{i=1}\mbox{Re}(\nabla
\psi^i,\psi^i)+4\sum\limits^4_{i=1}\mbox{Re}(\nabla\psi^i,\Pi\psi^i)\\
&&-2\sum\limits^4_{i=1}(\psi^i,\nu\cdot\nabla\nu\cdot\psi^i)
+2\sum\limits^4_{i=1}(\psi^i,\Pi(\nu\cdot\nabla\nu\cdot\psi^i)).
\end{array} $$
We define the function $\hat{p}$ by truncating $p$,
\[ \hat{p} \;=\; \min \left( p, (\frac{L^2}{4}-2)^2 \right)\;. \]
Then $\nabla \hat{p}(x)$ vanishes unless $p(x) \leq
(\frac{L^2}{4}-2)^2$. In this case, we have
\[ (\frac{L^2}{4}-2)^2 \;\geq\; p(x) \;\geq\; (\|\Pi_x\| - \|1\|)^2 \]
and thus $\|\Pi_x\| \leq \frac{L^2}{4}$. According to Lemma~4.1,
this implies that $|\psi^i| \leq L$ for all $i=1,\ldots,4$. We
conclude that
\[ \nabla \hat{p}(x) \not = 0 \spc \Longrightarrow \spc
|\psi^i(x)| \leq L. \] The last inequality allows us to estimate
$\nabla \hat{p}$ as follows, $$|\nabla \hat{p}|\leq
4L\sum\limits^4_{i=1}|\nabla\psi^i|+L^3\sum\limits^4_{i=1}|\nabla\psi^i|
+8L^2|h|+ 2L^4|h|\,\,$$ with $|h|^2=\sum\limits_{i,k}h^2_{ik}$.
Integration gives $$\begin{array}{lcl} \|\nabla \hat{p}\|^2_2
&\leq& 2 L^2(4+L^2)^2
(\sum\limits^4_{i=1}\|\nabla\psi^i\|^2_2+4L^2\|h\|^2_2) \\ &=&
8L^2(4+L^2)^2(4\pi E +\|h\|^2_2).
\end{array}
$$ The Sobolev inequality yields $$\|\hat{p}\|^2_6 \;\leq\;
\frac{48L^2(4+L^2)^2}{k^2} (4\pi E +\|h\|^2_2).$$ Hence
$\hat{p}(x)<\varepsilon$ except for $x\in U$, where the measure of
$U$ is bounded by $$\mu (U)^{1/3}\leq
\frac{48L^2(4+L^2)^2}{\varepsilon^2k^2}(4\pi E +\|h\|^2_2).$$
Clearly, on $M \setminus U$, also $p(x) < \varepsilon$. \QED
\end{Proof}

\section{Estimates of the Curvature Tensor}
We denote the curvature tensor of $\Sigma_N$ restricted to $M$ by
$\ov{R}^{\Sigma}_M:=i^\ast\ov{R}^{\Sigma}\in\Omega^2(M,\mbox{End}(\Sigma_{N}))$,
where $i$ is the natural inclusion $i:M\to N$. Recall that
$\ov{R}^\Sigma$ is related to the Riemannian curvature tensor
$\ov{R}$ by~(\ref{curv}). We denote the pull-back of $\ov{R}$ to
$M$ by $\ov{R}_M$ and define its norm by
$$|\ov{R}_M|^2=\sum\limits^3_{i,j=1}\sum\limits^3_{\alpha,\beta=0}(\ov{R}_{ij\alpha\beta})^2$$
We now derive a pointwise estimate for the curvature tensor in
terms of the system of Dirac spinors $\psi^i$.
\begin{Lemma}\label{lemma51}
$$(1-\|1-\Pi\|)\:|\ov{R}_M|^2 \;\leq\; 8
\sum\limits^4_{i=1}|\ov{\nabla}^2\psi^i|^2_{\otimes 2}\;. $$
\end{Lemma}
{\Proof} The identity
$\ov{R}^\Sigma_M(v,w)\psi=\ov{\nabla}^2\psi(v,w)-\ov{\nabla}^2\psi(w,v)$
immediately yields that $$|\ov{R}^\Sigma_M\psi|^2_{\otimes 2}
\;\leq\; 4|\ov{\nabla}^2\psi|^2_{\otimes 2} \;. $$ In order to
estimate the term on the left, we choose for given $x \in M$ an
orthonormal frame $(\nu,e_1,e_2,e_3)$ with $\ov{\nabla}e_i(x)=0$
and an orthonormal basis $(\phi_a)_{a=1,\ldots,4}$ of
$\Sigma_{N,x}$. Then for any linear map
$A\in\mbox{End}(\Sigma_{N,x})$,
$$\Tr(A\Pi)(x)=\sum\limits^4_{a=1}(\phi_a,A\:\Pi_x\phi_a)(x)=\sum_{i=1}^4
(\psi_i(x),A\psi_i(x))\;. $$ Thus
\begin{eqnarray}
\lefteqn{ \sum_{i=1}^4 |\ov{R}^\Sigma_M\psi^i|^2 \;=\;
\sum_{i=1}^{4} \sum_{j,k=1}^3
(\psi^i,\ov{R}^{\Sigma\ast}_M(e_j,e_k)\ov{R}^\Sigma_M(e_j,e_k)\psi^i)
} \\ &=& \sum\limits^3_{j,k=1}\Tr(\ov{R}^{\Sigma*}_M(e_j,e_k)\:
\ov{R}^\Sigma_M(e_j, e_k)\:\Pi) \\ &\geq& \sum\limits^3_{j,k=1}
\left(
\|\ov{R}^\Sigma_M(e_j,e_k)\|^2-\|\ov{R}^\Sigma_M(e_j,e_k)\:\ov{R}^{\Sigma\ast}_M(e_j,e_k)\|\:\|1-\Pi\|
\right). \label{cc}
\end{eqnarray}

Next we compute the appearing Hilbert-Schmidt norms.
\begin{eqnarray*}
\lefteqn{ \ov{R}^{\Sigma\ast}_M(e_i,e_j)\:\ov{R}^\Sigma_M(e_i,e_j)
}
\\ &=&\frac{1}{16} \left(-\sum\limits^3_{k,l,m,n=1}\ov{R}_{ijkl}
\ov{R}_{ijmn}e_k\cdot e_l\cdot e_m\cdot e_n \:+\: 2
\sum\limits^3_{k,l,m=1} \ov{R}_{ij0k}\ov{R}_{ijlm}\nu\cdot
e_k\cdot e_l\cdot e_m \right. \\ &&\left.\hspace*{1cm}
+2\sum\limits^3_{k,l,m=1} \ov{R}_{ijkl}\ov{R}_{ij0m}e_k\cdot
e_l\cdot\nu\cdot e_m
\:+\:4\sum\limits^3_{k,l=1}\ov{R}_{ijok}\ov{R}_{ij0m}\nu\cdot
e_k\cdot\nu\cdot e_m\right)\\
&=&\frac{1}{8}\left(\sum\limits^3_{k,l=1}\ov{R}_{ijkl}\ov{R}_{ijkl}
\:+\: 2\sum\limits^3_{k,l,m=1}\ov{R}_{ij0k}\ov{R}_{ijlm}\nu\cdot
e_k\cdot e_l\cdot e_m \:+\:
2\sum\limits^3_{k=1}\ov{R}_{ij0k}\ov{R}_{ij0k}\right)
\end{eqnarray*}
Since the trace of the second term vanishes, we conclude that
\begin{equation}\sum_{j,k=1}^3 \|\ov{R}^\Sigma_M(e_j, e_k)\|^2 \;=\;
\frac{1}{2}|\ov{R}_M|^2 \;.\label{cd}\end{equation} Moreover,
\begin{equation} \sum\limits^3_{j,k=1}
\|\ov{R}^\Sigma_M(e_j,e_k)\:\ov{R}^{\Sigma\ast}_M(e_j,e_k)\|
\;\leq\; \sum\limits^3_{j,k=1} \|\ov{R}^\Sigma_M(e_j,e_k) \|^2
\;\leq\; \frac{1}{2}|\ov{R}_M|^2 \;,\label{ce}\end{equation}
Substituting (\ref{cd}) and (\ref{ce}) into (\ref{cc}) gives the
result. \QED

\section{Integration by Parts}

In this section we derive an $L^2$ bound for the second derivative of
a solution of the boundary value problem~(\ref{br}).
The argument is similar to that given in~\cite{BF}.

\begin{Lemma}\label{lemma61} Suppose that $L$ satisfies the hypothesis
of Proposition~3.3. Then any solution $\Psi$ of the boundary value problem~(\ref{br})
satisfies the inequality
\begin{eqnarray*}
\lefteqn{ \int\limits_M \eta|\nabla^2\psi|^2 d\mu \;\leq\; c_1\:
\sup_{M}\left( |\Delta\eta|+|\nabla\eta||h|+\eta\:
(|R|+|h|^2+|\ov{\nabla}h| \right) E } \\ &&+c_2\:L \:\sup_{M}
\left( \eta\:(|\ov{\nabla}\ov{R}_M|+|h||\ov{R}_M|)\right) \sqrt{E}
\\ &&+c_3\:\frac{\sqrt{L+1}}{k} \:(\sup_{M}\eta)\;
\sqrt{\||h|^2+|\nabla h|\|_{6/5}} \;
\left\||\ov{\nabla}\ov{R}_M|+|h||\ov{R}_M \right\|_{5/12}
\;\sqrt{E}
\end{eqnarray*}
\end{Lemma}
\begin{Proof} A calculation similar to the one following (\ref{1et})
yields that
\begin{eqnarray}
|\ov{\nabla}^2\psi|^2
&=&\sum\limits_{j,i}\mbox{Re}(\ov{\nabla}^\ast_j\ov{\nabla}^2_{j,i}\psi,\ov{\nabla}_i\psi)
\label{6.1}\\ &&+\frac{1}{2}\Delta|\ov{\nabla}\psi|^2
\label{6.2}\\
&&-\frac{1}{2}\sum\limits_{j,i}(\ov{\nabla}_i\psi,\nu\cdot\ov{\nabla}^2_{j,j}\nu\cdot\ov{\nabla}_i\psi)\label{6.3}\\
&&-\sum\limits_{j,i}\mbox{Re}(\ov{\nabla}^2_{j,i}\psi,\nu\cdot\ov{\nabla}_j\nu\cdot\ov{\nabla}_i\psi)\label{6.4},
\end{eqnarray}
where $(e_1,\dots,e_n)$ is a smooth orthonormal frame on $M$.

In order to estimate the integral $\int\limits_M\eta|\ov{\nabla}^2\psi|d\mu$
with a positive test function $\eta \in C^2(M)$,we
consider the summands in the above equation separately.
Integrating by parts in~(\ref{6.2}) and using the decay properties of $\psi$, we obtain
\begin{equation*}
\frac{1}{2}\int\limits_M\eta\Delta|\ov{\nabla}\psi|^2d\mu=\frac{1}{2}\int\limits_M\Delta\eta|\ov{\nabla}\psi|^2d\mu\:\leq\:4
\pi \sup_{M}\Delta\eta\: E\:.
\end{equation*}

To estimate (\ref{6.3}) and (\ref{6.4}), we first calculate
\begin{eqnarray*}
\lefteqn{
\frac{1}{2}\ov{\nabla}_j(\ov{\nabla}_i\psi,\nu\cdot\ov{\nabla}_j\nu\cdot\ov{\nabla}_i\psi)
\;=\; \frac{1}{2}\ov{\nabla}_j \langle
\ov{\nabla}_i\psi,\ov{\nabla}_j\nu\cdot\ov{\nabla}_i\psi \rangle }
\\ &=&
\mbox{Re}(\ov{\nabla}^2_{j,i}\psi,\nu\cdot\ov{\nabla}_j\nu\cdot\ov{\nabla}_i\psi)
\;+\;
\frac{1}{2}(\ov{\nabla}_i\psi,\nu\cdot\ov{\nabla}^2_{j,j}\nu\cdot\ov{\nabla}_i\psi)
\;.
\end{eqnarray*}
Therefore, integration by parts gives
\begin{eqnarray*}
\lefteqn{ \left|\sum\limits_{i,j}\int\limits_M
  \eta\left(\frac{1}{2}\left(\ov{\nabla}_i\psi,\nu\cdot\ov{\nabla}^2_{j,j}
  \nu\cdot\ov{\nabla}_i\psi\right)+\mbox{Re}
\left(\ov{\nabla}^2_{j,i}\psi,\nu\cdot\ov{\nabla}_j\nu\cdot\ov{\nabla}_i\psi\right)\right)
d\mu \right| } \\ &\leq&
\frac{1}{2}\sum\limits_{i,j}\int\limits_M\left|(\partial_j\eta)\left(\ov{\nabla}_i\psi,\nu\cdot\ov{\nabla}_j
\nu\cdot\ov{\nabla}_i\psi\right)\right| d\mu \;\leq\; 4 \pi
\sup_{M}\left(|\nabla\eta||h| \right) E\;.
\end{eqnarray*}

It remains to control~(\ref{6.1}). Commuting the covariant derivatives,
we obtain, as in~\cite[eqns.~(31)-(35)]{BF},
\begin{eqnarray*}
\ov{\nabla}^\ast_j\ov{\nabla}^2_{j,i}\psi
&=&\ov{\nabla}_i(\ov{\nabla}^\ast\ov{\nabla}\psi)+\frac{1}{2}\ov{\nabla}^\ast_j(\ov{R}^\Sigma(e_j,e_i)\psi)+\\
&&+\mbox{Ric}(e_i,e_k)\ov{\nabla}_k\psi-\frac{1}{2}\ov{R}^\Sigma(e_j,e_j)\ov{\nabla}_j\psi
+\nu\cdot\ov{\nabla}^2_{i,j}\nu\cdot\ov{\nabla}_j\psi,
\end{eqnarray*}
where $\mbox{Ric}$ denotes the Ricci curvature of the hypersurface
$M\subset N$. If $\psi$ is a solution of~(\ref{br}), the first term can be simplified with
the Weitzenb{\"o}ck formula. Using the Gauss equation, we thus obtain
\[ \left(\ov{\nabla}^\ast_j\ov{\nabla}_{j,i}^2\psi,\ov{\nabla}_i\psi\right) \;\leq\;
\tilde{c}_1(|R|+|h|^2+|\ov{\nabla}h|)|\ov{\nabla}\psi|^2
\:+\: \tilde{c}_2\left(|\ov{\nabla}\ov{R}_i|+|h||\ov{R}_i|\right)\left(\psi,\ov{\nabla}_i\psi\right) \]
with suitable constants $\tilde{c}_1$ and $\tilde{c}_2$ which are
independent of the geometry. Now we choose $L>0$ as in Proposition~3.3 and calculate
\begin{eqnarray*}
\lefteqn{ \int\limits_M
\eta\,\mbox{Re}(\ov{\nabla}^\ast_j\ov{\nabla}^2_{j,i}\psi,\ov{\nabla}_i\psi)\:
d\mu \;\leq\; \tilde{c}_3\sup_{M}(\eta(|R|+|h|^2+|\ov{\nabla}h|)E
} \\ &&+\tilde{c}_2 \left( \int_{M\setminus\Omega_L} +
\int_{\Omega_L} \right) \eta
\left(|\ov{\nabla}^\ast\ov{R}_M|+|h|\:|\ov{R}_M| \right) |\psi|\:
|\ov{\nabla}\psi| \:d\mu \\ &\leq&
\tilde{c}_3\sup_{M}(\eta(|R|+|h|^2+|\ov{\nabla}h|)E \\ &&
+\tilde{c}_4 \:L\: \sup_{M} \left( \eta(|\ov{\nabla}\:\ov{R}_M)|
+|h|\:|\ov{R}_M|)\right) \sqrt{E} \\ &&+\tilde{c}_5
\int_{\Omega_L} \eta\: \left(
|\ov{\nabla}\:\ov{R}_M|+|h||\ov{R}_M| \right) \sqrt{|\psi|^2-1}\:
|\ov{\nabla}\psi|\:d\mu\;,
\end{eqnarray*}
where we have used the inequality $|\psi|\leq\sqrt{|\psi|^2-1}+1$.
In the last integral, we apply H{\"o}lder's inequality,
\begin{eqnarray*}
\lefteqn{ \int_{\Omega_L} \eta\: \left(
|\ov{\nabla}\:\ov{R}_M|+|h||\ov{R}_M| \right) \sqrt{|\psi|^2-1}\:
|\ov{\nabla}\psi|\:d\mu } \\ &\leq& \sup_M \eta\:
\left\||\ov{\nabla}\:\ov{R}_M|+|h||\ov{R}_M|
\right\|_{\frac{5}{12}}\: \left\||\psi|^2-1
\right\|^{\frac{1}{2}}_{L^6(\Omega_L)} \: \|\nabla \psi\|_2\;.
\end{eqnarray*}
Finally, the factor $\||\psi|^2-1\|_{L^6(\Omega_L)}$ is controlled by Proposition~3.3.
\end{Proof}
\QED

\noindent
{\em{Proof of Theorem~1.3.}}
For $L$ as in Proposition~3.3 and $\varepsilon=\frac{1}{2}$, we choose $U$ as in
Lemma~\ref{lemma42} to obtain
\[ \int\limits_{M\setminus U} \eta \:|\ov{R}_M|d\mu \;\leq\;
2\int\limits_{M\setminus U} \eta\:(1-\|1-\Pi\|)|\ov{R}_M|\:d\mu \;. \]
We now apply Lemma~5.1,
\[ \int\limits_{M\setminus U} \eta \:|\ov{R}_M|d\mu \;\leq\;
16\int\limits_M\eta\sum\limits^4_{i=1}|\ov{\nabla}^2\psi_i|^2_{\otimes2}\: d\mu \;.
\end{equation*}
Lemma~6.1 completes the proof. \QED

\newpage

Felix Finster

Naturwissenschaftliche Fakult\"{a}t I -- Mathematik

Universit\"{a}t Regensburg

D-93040 Regensburg

GERMANY

{\tt{felix.finster@mathematik.uni-regensburg.de}}

\medskip

Margarita Kraus

Naturwissenschaftliche Fakult\"{a}t I -- Mathematik

Universit\"{a}t Regensburg

D-93040 Regensburg

GERMANY

{\tt{margarita.kraus@mathematik.uni-regensburg.de}


\begin{thebibliography}{99}
\bibitem{ADM} R.\ Arnowitt, S.\ Deser, C.\ Misner, ``Energy and the
criteria for radiation in general relativity,'' {\em{Phys.\
Rev.}}\ 118, 1100 (1960)
\bibitem{Ba} R. Bartnik, "The Mass of an
Asymptotically Flat Manifold", {\em{Commun.\ Pure Appl.\ Math.}}\ XXXIX,
661-693 (1986)
\bibitem{Baum} H. Baum, ``Spin-Strukturen und Dirac-Operatoren {\"u}ber
pseudo-Riemannschen Mannigfaltigkeiten,'' {\em{Teubner Verlag Leipzig}} (1981)
\bibitem{BF} H.\ Bray, F.\ Finster, ``Curvature estimates and the positive
mass theorem,'' {\em{Comm.\ Anal.\ Geom.}}\ 10 (2002) 291-306
\bibitem{FK} F.\ Finster, I.\ Kath, ``Curvature estimates in asymptotically
flat manifolds of positive scalar curvature,'' {\em{Comm.\ Anal.\ Geom.}}\ 10
(2002) 1017-1031
\bibitem{Lami} H.-B. Lawson, M.-L. Michelsohn, {\it Spin Geometry},
Princeton University Press, Princeton 1989.
\bibitem{SY1} R.\ Schoen, S.-T.\ Yau, ``On the proof of the positive
mass conjecture in general relativity,'' {\em{Commun.\ Math.\ Phys.}}\ 65,
45-76 (1976)
\bibitem{PT} T.\ Parker, C.\ H.\ Taubes, ``On Witten's proof of the
positive energy theorem,'' {\em{Commun.\ Math.\ Phys.}}\ 84, 223-238 (1982)
\bibitem{W} E.\ Witten, ``A new proof of the positive energy
theorem,'' {\em{Commun.\ Math.\ Phys.}}\ 80, 381-402 (1981)
\end{thebibliography}
\end{document}
%%%Local Variables:
%%% mode: latex
%%% TeX-master: t
%%% End: